\documentclass[%
 aip,
 amsmath,amssymb,
 reprint,%
]{revtex4-1}

\usepackage{graphicx}
\usepackage{dcolumn}
\usepackage{bm}

\usepackage[utf8]{inputenc}
\usepackage[T1]{fontenc}
\usepackage{mathptmx}
\usepackage[mathscr]{euscript}
\usepackage{enumitem}

\begin{document}

\preprint{AIP/123-QED}

\title{On Koopman Mode Decomposition and Tensor Component Analysis}

\author{William T. Redman}
 \email{wredman@ucsb.edu}
 \affiliation{Interdepartmental Graduate Program in Dynamical Neuroscience, University of California Santa Barbara, California 93106, USA}
\

\date{\today}

\begin{abstract}
Koopman mode decomposition and tensor component analysis (also known as CANDECOMP/PARAFAC or canonical polyadic decomposition) are two popular approaches of decomposing high dimensional data sets into modes that capture the most relevant features and/or dynamics. Despite their similar goal, the two methods are largely used by different scientific communities and are formulated in distinct mathematical languages. We examine the two together and show that, under certain conditions on the data, the theoretical decomposition given by tensor component analysis is the \textit{same} as that given by Koopman mode decomposition. This provides a ``bridge'' with which the two communities should be able to more effectively communicate. Our work provides new possibilities for algorithmic approaches to Koopman mode decomposition and tensor component analysis, and offers a principled way in which to compare the two methods. Additionally, it builds upon a growing body of work showing that dynamical systems theory, and Koopman operator theory in particular, can be useful for problems that have historically made use of optimization theory. 
\end{abstract}

\maketitle

\begin{quotation}
Koopman operator theory (KOT) has emerged as a powerful and general framework with which to understand nonlinear dynamical systems in a data-driven manner. Despite its many strengths, there are numerous scientific domains, such as neuroscience and biology, where it remains an unfamiliar tool and therefore, is under employed. As part of this stems from KOT being mathematically formulated in a distinct way from other more commonly used techniques, such as principal component analysis (PCA), there is a need for bridging KOT with such methods. Here, we investigate how a major KOT analysis approach differs from tensor component analysis (TCA), an extension of PCA that has quickly become adopted by some of those same fields ``hesitant'' of KOT. We show that, in certain scenarios, TCA and KOT give the \textit{same} theoretical decomposition. This not only makes strides in establishing a bridge between TCA and KOT, but also provides a means with which to compare the strengths of the two methods and offers the possibility of new algorithmic approaches for KOT and TCA. 
\end{quotation}

\section{Introduction}
A central goal in all fields of science is to discover and interpret the underlying dynamics of observed data. Koopman operator theory (KOT) has emerged as a powerful, data-driven framework with which to do this\cite{koo31, koo32, mez05, bud12, mez19, bru21}. KOT lifts the perspective of dynamical systems from the possibly finite and nonlinear state space of the data, to the infinite and linear functional space of observables on the state space. A key finding over the past 15 years is that there are a variety of algorithms that can well approximate the Koopman operator, and hence the underlying dynamics, with sufficient, but reasonably finite, amounts of data. These algorithms include, among others, generalized Laplace analysis, the finite section method (also known as the Galerkin projection), Krylov subspace methods, and dynamic mode decomposition (see Mezić 2020 \cite{mez20} for a recent review). We note that all of these methods were originally developed, at least in part, independently of KOT. It has been the realization of their utility in computing relevant KOT quantities that has allowed for KOT to find success in an ever growing number of different scientific domains.

In many cases, researchers applying KOT are interested in the Koopman mode decomposition (KMD). The KMD provides a decomposition of the action of the Koopman operator in such a way that, for autonomous dynamical systems, the time evolution of the system is naturally given as the sum of a number of ``modes'', each with their own time dependencies (see Sec. II for more details). In some cases, the number of modes necessary for accurate prediction is small compared to the dimension of the state space, allowing for a low dimensional description. KMD has been successfully used to discover the underlying dynamics in many different dynamical systems such as fluids \cite{row09, mez13, arb17}, power grids \cite{sus16, kor18} and neural networks \cite{dog20, man20, tan20, nai21}. 

An active and fruitful area of research has been comparing the KMD modes to the modes obtained from various more commonly used methods, such as principle component analysis (PCA -- also referred to as proper orthogonal decomposition) and independent component analysis (ICA) \cite{mez05, row09, bru16, klu18a, lu20}. KMD and PCA have also been indirectly linked via their connections to the renormalization group  \cite{bra17, red20}. One common goal of this body of work has been to highlight to researchers in fields where PCA is a common tool that KMD can give similar and, in certain scenarios, superior mode representations of data. Despite this, to date KOT remains a niché and not widely adopted tool in some communities, such as neuroscience\cite{bru16, mar20} and biology\cite{bal20}. 

In recent years, a method that can be seen as a generalization to PCA, tensor component analysis\cite{car70, har70} (TCA -- also known as CADECOMP/PARAFAC or canonical polyadic decomposition), has gained in popularity. Part of this stems from the fact that TCA can be used on data that is naturally represented as a tensor, such as that obtained via repeated trial experiments. This means it does not require the trial averaging or ``flattening'' across dimensions that many matrix methods (e.g. PCA) do. Additionally, there exist a number of efficient and robust algorithms that have been developed to perform TCA and have been made freely available \cite{car70, har70, tensorlab, kol09, hon20, ram20}. Because of its intuitive connection to PCA, and because it does not require the computed components to be orthogonal to each other \cite{kru77} (unlike PCA), it has quickly been adopted by some of the same communities where KOT has not, such as neuroscience \cite{mar04, miw04, bec05, aca07, cic07, de_vos07, mor06, mor07, wil18, con19, xia20, zhu20a, zhu20b, xia21}. 

While there has been work extending KOT tools to tensors\cite{klu18b, gel19}, there has been little work comparing TCA to KMD\cite{lus19}. Here, by considering a data third-order tensor (i.e. a tensor with elements that are indexed by three independent labels), we show that there exists a correspondence between the computed modes of TCA and KMD. In particular, we show that, under certain conditions on the data, the decomposition obtained from applying TCA to data from an autonomous dynamical system is \textit{equivalent} to that obtained via KMD. We believe that this not only provides motivation as to why KMD can be used in fields that it is not yet popular in, but also highlights a new class of algorithms for computing the KMD (based on those from the TCA literature) and a new class of algorithms for performing TCA (based on those from the KMD literature). Because existing TCA and KMD algorithms rely on different approaches, understanding which ones work better under different conditions will be a fruitful future direction of study.

This paper is organized as follows. In Secs. II and III, we provide brief reviews of KMD and TCA respectively, with KMD being formulated so as to account for the tensor nature of the data, making its connection with TCA more apparent. We then proceed to layout the correspondence between the two methods in Sec. IV, and provide a clear description of when the two methods will (theoretically) provide exactly the same decompositions. We apply TCA and KMD to a simple numerical example to illustrate our claims in Sec. V. We end with Sec. VI where we discuss the limitations of our analysis, what our results allow us to say about how KMD and TCA compare, and what new questions emerge from our work. 

\section{Koopman Mode Decomposition}

The central object of interest in KOT is the Koopman operator, $\textbf{U}$, an infinite dimensional linear operator that describes the time evolution of observables (i.e. functions of the underlying state space variables) that live in the function space, $\mathcal{F}$. That is, after $t > 0$ amount of time, which can be continuous or discrete, the value of the observable $f \in \mathcal{F}$, which can be a scalar or a vector valued function, is given by
\begin{equation}
    \label{KO}
    \textbf{U}^t f(p) = f \left[ \textbf{T}^t(p) \right],
\end{equation}
where $\textbf{T}$ is the dynamical map evolving the system and $p$ is the initial condition or location in state space. For the remainder of the paper, it will be assumed that the state space being considered is of finite dimension and that $\mathcal{F}$ is the space of square-integrable functions.

The action of the Koopman operator on the observable $f$ can be decomposed as 
\begin{equation}
\label{KMD}
     \textbf{U} f(p) = \sum_{r = 0}^\infty \lambda_r \phi_r(p) \textbf{v}_r,
\end{equation}
where the $\phi_r$ are eigenfunctions of $\textbf{U}$, with $\lambda_r \in \mathbb{C}$ as their eigenvalues and $\textbf{v}_r$ as their Koopman modes \cite{mez05}. For systems with chaotic or shear dynamics, there is an additional term in Eq. \ref{KMD} arising from the continuous part of the spectrum of the Koopman operator \cite{mez05}. In general, not much is known about the contribution of this term \cite{bud12}. For the remainder of this paper, it will be assumed that the dynamical systems we are considering are such that the Koopman operator only has a point spectrum. 

Decomposing the action of the Koopman operator is powerful because, for a discrete dynamical system, the value of $f$ at time step $k \in \mathbb{Z}^+$ is given simply by
\begin{equation}
    \label{KMD time}
     f\left[ \textbf{T}^k(p) \right]  = \textbf{U}^k f(p) = \sum_{r = 0}^\infty \lambda_r^k \phi_r(p) \textbf{v}_r.
\end{equation}
When the dynamical system is continuous, the value of $f$ at time $t \in \mathbb{R}^+$ is given by
\begin{equation}
    \label{KMD cont}
     f\left[ \textbf{T}^t(p) \right] = \textbf{U}^t f(p) = \sum_{r = 0}^\infty \exp(\lambda_r t) \phi_r(p) \textbf{v}_r.
\end{equation}
Note here that the convention changes from $\lambda_r^k$ to $\exp(\lambda_r t)$.

From Eqs. \ref{KMD time} and \ref{KMD cont}, we see that the dynamics of the system in the directions $\textbf{v}_r$, scaled by $\phi_r(p)$, are given by the magnitude of the corresponding $\lambda_r$. Assuming that $|\lambda_r| \leq 1$ for all $r$, finding the long time behavior of $f$ amounts to considering only the $\phi_r(p)\textbf{v}_r$, whose $|\lambda_r| \approx 1$ (in the case of a discrete dynamical system).  We note that the eigenfunctions, $\phi_r$, are the only parts of Eqs. \ref{KMD time} and \ref{KMD cont} that retain any ``memory'' of the initial condition \cite{bud12}. 
 
While the number of triplets $(\textbf{v}_r, \lambda_r, \phi_r)$ needed to fully capture the action of $\textbf{U}$ is, in principle, infinite, in many applied settings it has been found that a finite number, $R \in \mathbb{Z}^+$, of them allows for a good approximation \cite{bud12}. That is, in the case of a discrete dynamical system,
\begin{equation}
    \label{KMD finite}
    \textbf{U}^k f(p) \approx \sum_{r = 0}^{R - 1} \lambda_r^k \phi_r(p) \textbf{v}_r.
\end{equation}
In this paper, we consider Eq. \ref{KMD finite} to be the definition of KMD. 

\subsection{Dynamic Mode Decomposition}
One popular way of computing the triplets $(\textbf{v}_r, \lambda_r, \phi_r)$ of Eq. \ref{KMD finite} is dynamic mode decomposition (DMD) \cite{sch10, row09, tu14}. While a number of different variants of DMD have been developed\cite{sch10, che12, tu14, wil15, arb17}, we will consider Exact DMD, which was first proposed in Tu et al. 2014\cite{tu14}. 

Let $\textbf{X} = [\textbf{z}_0, \textbf{T} \textbf{z}_0, ..., \textbf{T}^{m-1} \textbf{z}_0] = [\textbf{z}_0, \hspace{1mm} \textbf{z}_1,..., \textbf{z}_{m-1}] \in \mathbb{R}^{n \times m}$ be $m$ snapshots of the state space column vector $\textbf{z} \in \mathbb{R}^n$. Here $\textbf{T}$ is the dynamical map evolving the system. Let $\textbf{Y} = [\textbf{z}_1, \hspace{1mm} \textbf{z}_2, ..., \textbf{z}_{m}]$. Exact DMD does not require $\textbf{X}$ and $\textbf{Y}$ to be built from sequential time-series data, but we assume it because it will make the connection to TCA more clear.

Exact DMD considers the operator $\textbf{A}$, defined as 
\begin{equation}
    \label{Exact DMD}
    \textbf{A} = \textbf{Y} \textbf{X}^+,
\end{equation}
where $^+$ is the Moore-Penrose pseudoinverse. Eq. \ref{Exact DMD} is the least-squares solution to $\textbf{A} \textbf{X} =\textbf{Y}$, and hence, satisfies
\begin{equation}
    \label{Least Squares}
   \min_{\textbf{A}} \Big|\Big|\textbf{Y} - \textbf{A} \textbf{X} \Big|\Big|_F,
\end{equation}
where $|| \cdot ||_F$ is the Frobenius norm. If $\textbf{Yc} = \textbf{0}$ for all $\textbf{c} \in \mathbb{R}^{n}$ such that $\textbf{Xc} = \textbf{0}$, then $\textbf{A}\textbf{X} = \textbf{Y}$ exactly \cite{tu14}. This property is called linear consistency.  

In practice, the DMD modes are found by computing the (reduced) singular value decomposition (SVD) of $\textbf{X}$, $\textbf{X} = \textbf{Q} \bm{\Sigma} \textbf{V}^*$, and setting $\tilde{\textbf{A}}=\textbf{Q}^* \textbf{Y} \textbf{V} \bm{\Sigma}^{-1}$ (see Algorithm 2 and Appendix 1 of Tu et al. 2014\cite{tu14} for more details on how the modes are computed from $\tilde{\textbf{A}}$). Importantly, when $\textbf{X}$ and $\textbf{Y}$ are linearly consistent and $\textbf{A}$ has a full set of eigenvectors, the DMD modes are equivalent to the Koopman modes\cite{tu14}. 

Let $\textbf{A}$ have a full set of eigenvectors $\textbf{v}_r$ with eigenvalues $\lambda_r$. We can then write $\textbf{x}_0$, the first column of $\textbf{X}$, as
\begin{equation}
    \label{First column of X}
    \textbf{x}_0 = \sum_{r = 0}^{R - 1} \phi_r(\textbf{z}_0) \textbf{v}_r,
\end{equation}
for some choice of constants $\phi_r(\textbf{z}_0)$, which correspond to the Koopman eigenfunctions. Here $R$ is the rank of $\tilde{\textbf{A}}$. If $\textbf{X}$ and $\textbf{Y}$ are linearly consistent, then -- given the relationship between the columns of $\textbf{X}$ -- any column of $\textbf{X}$ can be written as
\begin{equation}
    \label{Other columns of X}
    \textbf{x}_k = \textbf{A}^k\textbf{x}_0 = \sum_{r = 0}^{R - 1} \lambda_r^k\phi_r(\textbf{z}_0) \textbf{v}_r,
\end{equation}
for $k < m - 1$. Similarly, any column of $\textbf{A}\textbf{X}$ can be written as
\begin{equation}
    \label{Columns of UX} \left[ \textbf{A}\textbf{X}\right]_k = \textbf{A}\textbf{x}_k = \sum_{r = 0}^{R - 1} \lambda_r^{k+1}\phi_r(\textbf{z}_0) \textbf{v}_r.
\end{equation}
Defining the matrix 
\begin{equation*}
\textbf{S} =
    \begin{bmatrix} 
    \lambda_0 & \lambda_1 & \cdot \cdot \cdot & \lambda_{R-1}\\
    \lambda_0^2 & \lambda_1^2 & \dots & \lambda_{R-1}^2\\
    \vdots & \vdots & \ddots & \vdots \\
    \lambda_{0}^m & \lambda_{1}^m & \dots & \lambda_{R-1}^{m}\\
    \end{bmatrix},
\end{equation*}
we can write
\begin{equation}
\label{dmd single}
    \textbf{Y} - \textbf{A} \textbf{X} = \textbf{Y} - \sum_{r = 0}^{R-1}  \left[\phi_r(\textbf{z}_0)\textbf{v}_r\right] \otimes \textbf{s}_r = \textbf{0},
\end{equation}
where $\textbf{s}_r$ denotes the $r^{\text{th}}$ column of $\textbf{S}$, $\otimes$ is the vector outer product, and the second equality is due to the assumed linear consistency of the data. 

In some cases, DMD is estimated using only a single time series. In such a case, the computed eigenfunctions $\phi_r(\textbf{z}_0)$ are computed for only a single $\textbf{z}_0$. To gain more information on the $\phi_r$, multiple experiments, each with a different initial condition, can be performed. Alternatively, multiple sensors or sources of data (e.g. physical locations in a fluids experiment) can be considered as different initial conditions. These may not be feasible for a given experimental design, especially when analyzing previously recorded data. In general, there exists a way in which to discover more values of the eigenfunctions. For deterministic systems, time delaying \cite{tak81, kam20} a given time series (e.g. considering the location of the system at $t_1$ to be the initial condition, instead of that at $t_0$, and so on) provides more initial conditions where the eigenfunctions can be evaluated at. This is, of course, limited by the fact that only the points that fall along the single trajectory that the data came from can be used.

When data from $q$ different initial conditions (obtained via whichever of the methods listed above), $\{\textbf{z}^{(0)}_0, ..., \textbf{z}_0^{(q-1)}\}$, is considered, the data can be represented by third-order tensors, $\pmb{\mathscr{X}}$ and $\pmb{\mathscr{Y}}$. To do this, data matrices, as we had before, can be ``stacked'' along a third dimension. The $j^{\text{th}}$ frontal slice of $\pmb{\mathscr{X}}$ (i.e. the matrix made by fixing the third index of the tensor to be $j$) is therefore given by the data matrix $\textbf{X}_j$. Here $\textbf{X}_j$ corresponds to the data collected with initial condition $\textbf{z}_0^{(j)}$. As noted earlier, the dependence on $\textbf{z}_0^{(j)}$ in Eq. \ref{KMD finite} is only in the eigenfunctions $\phi_r$. If we define ${\bm \varphi}_r^T = \big[\phi_r(\textbf{z}_0^{(0)}), ..., \phi_r(\textbf{z}_0^{(q - 1)})\big]$, then we have an equation analogous to Eq. \ref{dmd single},
\begin{equation}
    \label{dmd multiple} 
    \pmb{\mathscr{Y}} - \sum_{r = 0}^{R-1} \textbf{v}_r  \otimes \textbf{s}_r \otimes {\bm \varphi}_{r} = \textbf{0},
\end{equation}
again assuming that all corresponding pairs of $\textbf{X}_j$ and $\textbf{Y}_j$ are linearly consistent. 

While all of this has been shown for $\pmb{\mathscr{X}}$ and $\pmb{\mathscr{Y}}$ being made up of snapshots of the state space vector $\textbf{z}$, the same holds true for the case when the data tensor is instead composed of other, non-identity, observables \cite{row09, tu14}. There are many cases where this has been found to give better decompositions for the right choices of observables. 

\section{Tensor Component Analysis}
Similar to KOT, TCA was first developed in the early $20^{\text{th}}$ century \cite{hit27, hit28} and has seen a recent resurgence of interest. It has been successfully applied to a number of fields including neuroscience \cite{mar04, miw04, bec05, aca07, cic07, de_vos07, mor06, mor07, wil18, con19, xia20, zhu20a, zhu20b, xia21} and signal processing \cite{mut05, de_lat07, de_lat08}, as well as the original fields it was developed in, psychometrics \cite{car70, har70, tuc66} and chemometrics  \cite{app81, leu92, hen93, smi94, kie98, Wu98, bro99, and03, smi04}. Given that it was independently arrived upon by several different researchers at different times in the context of different fields, the TCA literature is somewhat confusing. TCA is also referred to as CANDECOMP (for \underline{can}onical \underline{decomp}osition), PARAFAC (for \underline{para}llel \underline{fac}torization), canonical polyadic decomposition, and CP decomposition.  For a comprehensive review of the literature on tensor decompositions, see Kolda and Bader 2009\cite{kol09}. 

Given a tensor $\pmb{\mathscr{Y}}$, the central objective of TCA is to find vectors that well reconstruct it. While $\pmb{\mathscr{Y}}$ may come from the data of a dynamical system, it does not have to. For example, $\pmb{\mathscr{Y}}$ can be populated by the pixel values of a set of images. The goal of applying TCA, in such a case, could be to find features common in different subsets of the images. In the case of $\pmb{\mathscr{Y}}$ being a third-order tensor of size $I_1 \times I_2 \times I_3$, TCA finds matrices $\textbf{A} \in \mathbb{C}^{I_1 \times R}, \hspace{0.5mm} \textbf{B} \in \mathbb{C}^{I_2 \times R}, \hspace{0.5mm}  \textbf{C} \in \mathbb{C}^{I_3 \times R}$, such that
\begin{equation}
\label{TCA eq}
    \pmb{\mathscr{Y}} \approx  \sum_{r = 0}^{R - 1} \textbf{a}_r \otimes \textbf{b}_r \otimes \textbf{c}_r,
\end{equation}
where $\textbf{a}_r$, $\textbf{b}_r$, and $\textbf{c}_r$ are, respectively, the $r^{\text{th}}$ column vectors of $\textbf{A}$, $\textbf{B}$, and $\textbf{C}$, $\otimes$ is the vector outer product, and $R \in \mathbb{Z}^+$. The TCA modes are therefore given by the triplets $(\textbf{a}_r, \textbf{b}_r, \textbf{c}_r)$. $\textbf{A}, \textbf{B}$, and $\textbf{C}$ are found by the minimization problem
\begin{equation}
    \label{TCA min}
    \min_{\textbf{A},  \textbf{B}, \textbf{C}} \Big|\Big| \pmb{\mathscr{Y}} - \sum_{r = 0}^{R - 1} \textbf{a}_r \otimes \textbf{b}_r  \otimes \textbf{c}_r \Big|\Big|,
\end{equation}
where $|| \cdot ||$ is analogous to the Frobenius norm for matrices \cite{kol09},
\begin{equation*}
    \Big|\Big| \pmb{\mathscr{Y}} \Big|\Big| = \sqrt{\sum_{i_1 = 0}^{I_1 - 1} \hspace{0.75mm} \sum_{i_2 = 0}^{I_2 - 1} \sum_{i_3 = 0}^{I_3 - 1} x^2_{i_1 i_2 i_3}} \hspace{1mm}.
\end{equation*}
The rank of $\pmb{\mathscr{Y}}$, denoted $R^*(\pmb{\mathscr{Y}}) =$ rank($\pmb{\mathscr{Y}}$), is defined as the minimum number of triplets needed for the solution to Eq. \ref{TCA min} to be zero \cite{kol09}. The decomposition obtained using $R^*(\pmb{\mathscr{Y}})$ is called the \textit{rank decomposition} of $\pmb{\mathscr{Y}}$. While this is, in principle, the most natural definition for the TCA representation of $\pmb{\mathscr{Y}}$, computing $R^*(\pmb{\mathscr{Y}})$ can be a highly non-trivial problem\cite{kol09}. While there are some results on the maximal and most likely ranks to occur given the dimensionality of $\pmb{\mathscr{Y}}$, they are limited \cite{kru89, ber99, ber00}. Therefore, in practical settings, many different values of $R$ are tried, and the smallest one that satisfies some criteria (e.g. sufficiently small reconstruction error) is chosen. For sake of clarity, we refer to the decompositions that satisfy Eq. \ref{TCA min} for a given $R > R^*(\pmb{\mathscr{Y}})$ as \textit{R--optimal TCA decompositions} of $\pmb{\mathscr{Y}}$.

\textbf{Lemma 1} \textit{R--optimal TCA decompositions reconstruct $\pmb{\mathscr{Y}}$ exactly.} 

\textit{Proof.} Let  $\textbf{A}, \textbf{B}, \textbf{C}$ be the matrices that correspond to the rank decomposition of $\pmb{\mathscr{Y}}$, and $\Delta R = R - R^*(\pmb{\mathscr{Y}})$. Let $\{\alpha_0,..., \alpha_{\Delta R} \}$ be such that $\alpha_i \in \mathbb{R}/\{0\}$, $\alpha_i \neq \alpha_j$, and $\sum_i \alpha_i = 1$. Then the matrices $\textbf{A}' = [\textbf{a}_0, ...,\textbf{a}_{R^*-2}, \alpha_0 \textbf{a}_{R^* - 1}, ..., \alpha_{\Delta R}\textbf{a}_{R^* - 1}]$, $\textbf{B}' = [\textbf{b}_0, ..., \textbf{b}_{R^* - 2}, \textbf{b}_{R^* - 1}, ..., \textbf{b}_{R^* - 1}]$, and $\textbf{C}' = [\textbf{c}_0, ..., \textbf{c}_{R^* - 2},\textbf{c}_{R^* - 1} ..., \textbf{c}_{R^* - 1}]$ make up a decomposition that has zero error in reconstructing $\pmb{\mathscr{Y}}$: 
\begin{equation*}
\begin{split}
    \sum_{r = 0}^{R-1} \hspace{2.5mm} & \textbf{a}'_r \otimes \textbf{b}'_r \otimes \textbf{c}'_r = \\
     \sum_{r = 0}^{R^*(\pmb{\mathscr{Y}})-2} & \textbf{a}_r \otimes \textbf{b}_r \otimes \textbf{c}_r + \sum_{i = 0}^{\Delta R} \alpha_i \textbf{a}_{R^*(\pmb{\mathscr{Y}}) - 1} \otimes \textbf{b}_{R^*(\pmb{\mathscr{Y}}) - 1} \otimes \textbf{c}_{R^*(\pmb{\mathscr{Y}}) - 1} =\\
    \sum_{r = 0}^{R^*(\pmb{\mathscr{Y}})-2} & \textbf{a}_r \otimes \textbf{b}_r \otimes \textbf{c}_r + \left[\textbf{a}_{R^*(\pmb{\mathscr{Y}}) - 1} \otimes \textbf{b}_{R^*(\pmb{\mathscr{Y}}) - 1} \otimes \textbf{c}_{R^*(\pmb{\mathscr{Y}}) - 1}\right] =\\
    \sum_{r = 0}^{R^*(\pmb{\mathscr{Y}})-1} & \textbf{a}_r \otimes \textbf{b}_r \otimes \textbf{c}_r = \pmb{\mathscr{Y}} \\
  \end{split}
\end{equation*}
Given that at least one R--optimal decomposition reconstructs $\pmb{\mathscr{Y}}$ exactly, from Eq. \ref{TCA min} we have that all  R--optimal decompositions must. $\square$

From the proof of Lemma 1, it is straightforward to see that there are many possible R--optimal decompositions and thus, R--optimal decompositions are not unique. However, rank decompositions can be unique. This is a particularly appealing feature of TCA, as it is not generally true of matrix decomposition methods (e.g. PCA). In particular, it has been proven that if
\begin{equation}
    \label{Kruskal uniqueness}
    \text{rank}(\textbf{A}) + \text{rank}(\textbf{B}) + \text{rank}(\textbf{C}) \geq 2\hspace{0.5mm} R^*(\pmb{\mathscr{Y}}) + 2,
\end{equation}
then the decomposition is unique, up to rescaling and permutation of the $(\textbf{a}_r, \textbf{b}_r, \textbf{c}_r)$\cite{kru77}. The non-uniqueness to scaling emerges because, for $\alpha \in \mathbb{C} / \{0\}$, $(\textbf{a}_r / \alpha, \textbf{b}_r, \alpha\textbf{c}_r)$ also satisfies Eq. \ref{TCA min}. Similarly, if the columns of $\textbf{A}, \textbf{B},$ and $\textbf{C}$ are shuffled in the same way, the new matrices $\tilde{\textbf{A}}, \tilde{\textbf{B}},$ and $\tilde{\textbf{C}}$ satisfy Eq. \ref{TCA min}.

As has been noted in the literature, there is no perfect way of performing TCA \cite{kol09}. In practice, various optimization based methods, such as alternating and nonlinear least squares, are used for solving Eq. \ref{TCA min}. These methods can take many iterations to converge and are known to not always be guaranteed to converge to a global minimum. Therefore, the development of methods to better compute TCA decompositions is an area of active research. More details on these and other algorithms, as well as their practical considerations, can be found in Kolda and Bader 2009 \cite{kol09} and Hong, Kolda, and Bader 2020 \cite{hon20} (which includes discussion on the use of more general objective functions). 

We note that, in view of Eq. \ref{TCA min}, TCA has usually been seen as a tool to represent $\pmb{\mathscr{Y}}$. That is, the modes have been used to provide insight into existing data, as opposed to predicting the future state(s) of the system. 

\section{Correspondence between TCA and KMD}
We now proceed to show how TCA and KMD are related. Let the third-order tensors $\pmb{\mathscr{X}}$ and $\pmb{\mathscr{Y}}$ be related to each other as in Sec. IIA and be made up of linearly consistent frontal slices (now referred to as data matrices). Let the approximated Koopman operator have rank $R$ and a full set of eigenvectors. Let $R^*(\pmb{\mathscr{Y}}) = \text{rank}(\pmb{\mathscr{Y}})$, where $R \geq R^*(\pmb{\mathscr{Y}})$. From Eq. \ref{dmd multiple} and Lemma 1, we have that the Exact DMD and TCA modes, $(\textbf{v}_r, \textbf{s}_r, \bm{\varphi}_r)$ and $(\textbf{a}_r, \textbf{b}_r, \textbf{c}_r)$ respectively, are such that   
\begin{equation}
    \label{general tca kmd comp}
    \begin{split}
        \pmb{\mathscr{Y}} - \sum_{r = 0}^{R-1}\textbf{v}_r \otimes \textbf{s}_r \otimes {\bm \varphi}_r  &= \textbf{0}\\
         \pmb{\mathscr{Y}} - \sum_{r = 0}^{R-1} \textbf{a}_r \otimes \textbf{b}_r \otimes \textbf{c}_r &= \textbf{0}.
    \end{split}
\end{equation}

\textbf{Lemma 2} \textit{Let $\pmb{\mathscr{X}}$ and $\pmb{\mathscr{Y}}$ be made up of linearly consistent data matrices. When $R > R^*(\pmb{\mathscr{Y}})$, the DMD modes are an R--optimal TCA decomposition.}

\textit{Proof.} This follows immediately from the fact that Eq. \ref{dmd multiple} implies that the DMD modes satisfy the TCA minimization problem (Eq. \ref{TCA min}). $\square$

Based on Eq. \ref{general tca kmd comp} and Lemma 2, we suggest the following correspondence:
\begin{equation}
\label{TCA KMD corr}
\begin{split}
  \textbf{v}_r &\Longleftrightarrow \textbf{a}_r  \\
  \textbf{s}_r&\Longleftrightarrow \textbf{b}_r \\
   {\bm \varphi}_r   &\Longleftrightarrow \textbf{c}_r  
\end{split}
\end{equation}

When, if ever, will this correspondence be exact? That is, when will the computed DMD triplets $(\textbf{v}_r, \textbf{s}_r, {\bm \varphi}_r)$, in theory, be equal to the computed TCA triplets $(\textbf{a}_r, \textbf{b}_r, \textbf{c}_r)$? 

\textbf{Lemma 3} \textit{Let $\pmb{\mathscr{X}}$ and $\pmb{\mathscr{Y}}$ be made up of linearly consistent data matrices. Let $\textbf{V} = [\textbf{v}_1, ..., \textbf{v}_R]$, $\textbf{S} = [\textbf{s}_1, ..., \textbf{s}_R]$, and $\bm{\Phi} = [\bm{\varphi}_1, ...,\bm{\varphi}_R]$. When $R = R^*(\pmb{\mathscr{Y}})$ and $\text{rank}(\bm{\Phi}) + \text{rank}(\textbf{V}) + \text{rank}(\textbf{S}) \geq 2R + 2$, then the correspondence of Eq. \ref{TCA KMD corr} is exact, up to scaling and permutation of labels.}

\textit{Proof.} This follows immediately from the fact that Eq. \ref{dmd multiple} implies that the DMD modes satisfy the TCA minimization problem (Eq. \ref{TCA min}) and that the TCA decomposition is unique up to scaling and permutation when $R = R^*(\pmb{\mathscr{Y}})$ and Eq. 15 is satisfied \cite{kru77}, which have been assumed to be true. $\square$

Therefore, there exists a restricted case (but a case nonetheless) where the TCA and DMD modes are equivalent! Lemma 3 additionally tells us that this equivalence will be met when the data comes from distinct sources that have single exponential growth/decay and/or oscillatory dynamics. This is because, for autonomous discrete dynamical systems, the time evolution of the $r^{\text{th}}$ KMD mode is described by integer powers of a single complex number, $\lambda_r$. While, \textit{a priori}, the TCA modes have no assumed form, when Eq. 17 is exact, their time dependence must match those of the DMD modes.

Note that if there exists multiple TCA decompositions [due to either $R^* > R(\pmb{\mathscr{Y}})$ or Eq. \ref{Kruskal uniqueness} not being satisfied], Lemma 2 tells us that one of them is equivalent to the DMD modes. In such a case, it would be well motivated to choose, as convention, the DMD  modes as \textit{the} TCA decomposition, given DMD's connection to the dynamics underlying the data. This decomposition will also have modes with only exponential and/or oscillatory time dependencies. 


\section{Numerical Example}

To illustrate the correspondence developed in Sec. IV, and to show that standard numerical implementations of TCA can indeed give decompositions that are equivalent to Exact DMD, we applied both methods to ``synthetic'' data. Following the common practice of testing new tensor decomposition methods on random tensors was not feasible here, because of the assumptions we have had to introduce. However, we have taken that approach in spirit by generating data that has both the necessary dynamical systems structure and elements that are chosen randomly. The tensors were constructed as follows. 

First, vectors $\{ \bar{\textbf{v}}_0,...,\bar{\textbf{v}}_{R - 1}\}$, with  $\bar{\textbf{v}}_r \in (0, 1)^{n}$ sampled uniformly, were generated and then made to be orthonormal. Next, exponential functions $\{\exp(2\pi i \alpha_0/10 + \beta_0/10), ...,\exp(2\pi i \alpha_{R-1}/10 + \beta_{R - 1}/10) \}$, with $\alpha_r, \beta_r \in (-1,1)$  sampled uniformly, were generated. Vectors $\bar{\textbf{s}}_r^T = [\exp(2\pi i \alpha_r/10 + \beta_r/10),...,\exp(2\pi i N \alpha_r/10 + \beta_r/10)]$, for a chosen $N \in \mathbb{Z}^+$, were created. Finally, random polynomial functions $\{\bar{\phi}_0,...,\bar{\phi}_{R-1} \}$, with $\bar{\phi}_r(x) = x^{\gamma_r}$ and $\gamma_r \in (0, K)$ sampled uniformly, were generated. Vectors $\bar{\bm{\varphi}}_r^T = [\bar{\phi}_r(0.1),...,\bar{\phi}_r(1)]$ were then created. 

Letting $\bar{\textbf{s}}'_r$ be the time shifted version of $\bar{\textbf{s}}_r$, i.e. $\bar{\textbf{s}}^{'T}_r = \exp(2\pi i \alpha_r/10 + \beta_r/10) \bar{\textbf{s}}_r^T$, the tensors
\begin{equation*}
\begin{split}
    \pmb{\mathscr{X}} = \sum_{r = 0}^{R-1}\bar{\textbf{v}}_r \otimes \bar{\textbf{s}}_r \otimes \bar{{\bm \varphi}}_r \\
     \pmb{\mathscr{Y}} = \sum_{r = 0}^{R-1}\bar{\textbf{v}}_r \otimes \bar{\textbf{s}}'_r \otimes \bar{{\bm \varphi}}_r
    \end{split}
\end{equation*}
were constructed. We note here that this data can be interpreted as coming from a linear dynamical system, where the $\bar{\textbf{v}}_i$ are the eigenvectors of the dynamical map, the $\bar{\textbf{s}}_i$ are the time evolution of each mode, and the $\bar{\phi}_i$ are the dependence of the amplitude of each mode on initial condition. Note the bars have been used to distinguish these sources from the computed Exact DMD modes (although we expect, and do indeed find, that they should be the same).  

The TCA modes were computed using Tensorlab 3.0\cite{tensorlab}, a MATLAB package that implements various possible TCA decomposition algorithms. We used the nonlinear least squares based method. The Exact DMD modes were computed using Algorithm 2 and Appendix 1 from Tu et al. 2014\cite{tu14}. The code that generated these examples has been made freely available online \cite{wtr_git}.

\begin{figure}
    \centering
    \includegraphics[width = 0.49\textwidth]{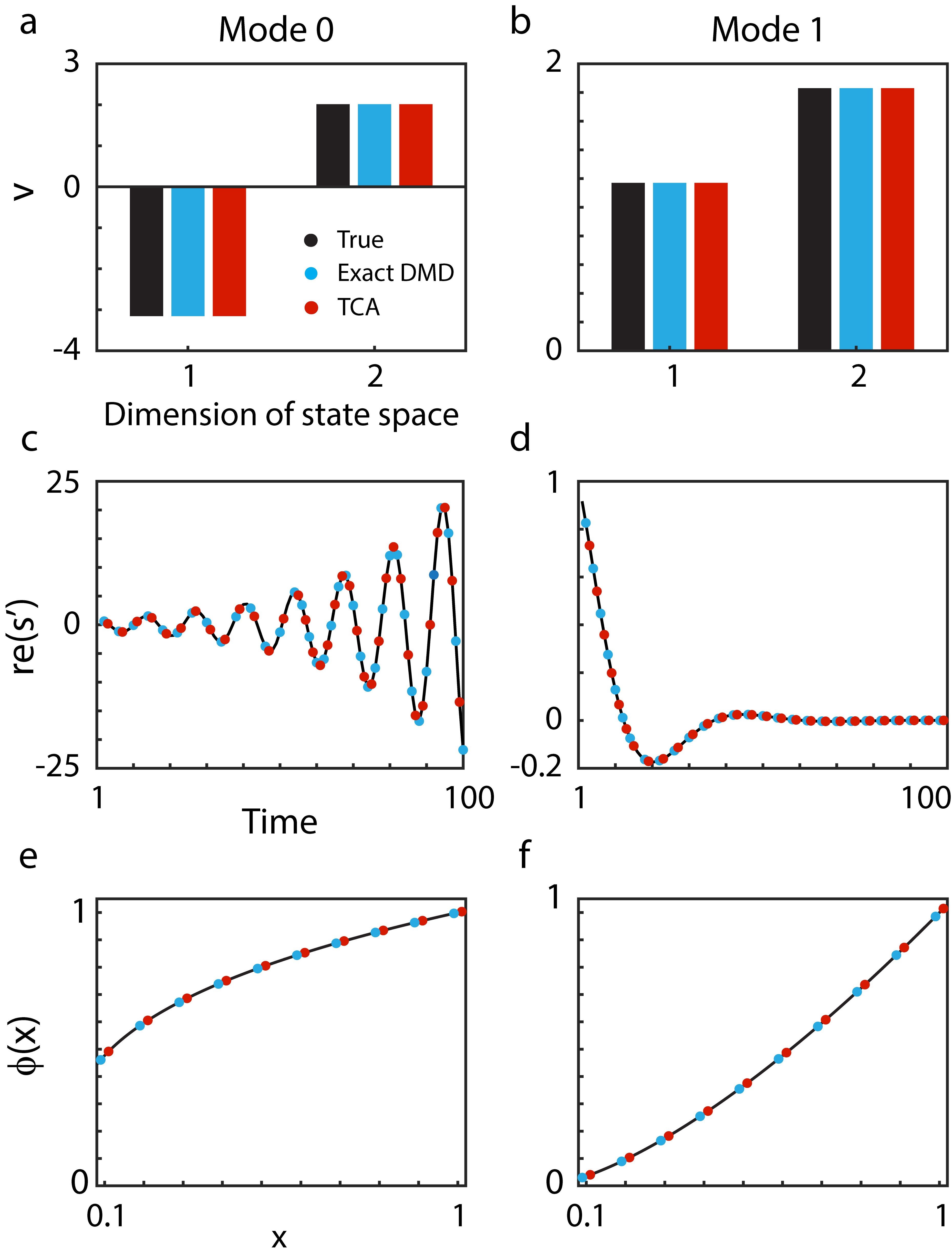}
    \caption{\textbf{Numerical example of TCA and DMD modes being equivalent.} (a) The two components of the vector $\textbf{v}_0$ as computed by Exact DMD (light blue) and TCA (red), along with their true values (black). The Exact DMD values are the mean of those computed from the data matrices corresponding to the $10$ different initial conditions. (c) The real part of each mode's time evolution $\textbf{s}'_0$, as computed by the two methods, along with their true value. The Exact DMD values are again the mean over the $10$ different computed values. Points are offset from each other and ``sparsefied'' for display purposes (they would otherwise lie on top of each other). (e) The dependence on initial condition $\phi_0(x)$, as computed by the two methods, along with their true values. The points are again separated for purpose of display. (b), (d), (f) are the same as (a), (c), (e), respectively, for the other mode. }
    \label{r2_lin_sys}
\end{figure}

The computed modes from one instance of the randomly generated data using $R = 2$, $n = 2$, $N = 100$, and $K = 3$ is shown in Fig. \ref{r2_lin_sys}. Here $\textbf{v}_0^T = [-3.1530, 2.0160]$, $\textbf{v}_1^T = [1.1704, 1.8305]$, $\alpha_0 = 0.7418$, $\alpha_1 = -0.2043$, $\beta_0 = 0.322$, $\beta_1 = -0.7883$, $\gamma_0 = 0.3217$, and $\gamma_1 = 1.4534$. For both modes, the TCA and Exact DMD triplets were nearly identical (they have been separated from each other for clarity of display), and lie on top of the true sources. Thus, both methods correctly reveal that the first mode is less sensitive to initial condition than the second (sublinear vs. superlinear), and that the first mode grows in time, whereas the second decays. 

We generated 100 data tensors (again using $R = 2$, $n = 2$, $N = 100$, and $K = 3$) and computed the mean normed difference between the Exact DMD and the TCA modes, normalizing by the norm of the TCA modes. The median of these 100 values was $2.43\times10^{-12}$. As this suggests, for many of the 100 tensors, the two methods produced very similar results. However, we did find a few cases where the difference between the two was quite large (for these cases, the difference between the Exact DMD modes and the true values remained very small). The extent of this appeared to be related to how ``long'' the experiment was, as increasing $N$ led to more discrepancies and decreasing $N$ led to more overlap. We do not believe this was (primarily) caused by a change in $R^*$, as the error associated with using only one mode remained far larger (many orders of magnitude) than the error associated with using two. A reason for these results may be that, by increasing $N$ (and thus, the size of the matrix $\textbf{B}$), more local minima emerge that the optimization algorithm could get trapped in. 

\section{Discussion}
In this paper, we examined tensor component analysis (TCA -- also known as  CANDECOMP/PARAFAC or canonical polyadic decomposition) in relation to Koopman mode decomposition (KMD). This was motivated by the fact that both methods have become popular ways to discover, in an unsupervised manner, the relevant features and/or dynamics of a given dynamical system. Despite their joint aim, the two methods have largely occupied disjoint scientific realms. Therefore, it became our goal to examine the two methods together and see what, if any, connections existed between them in an effort to ``bridge'' the different communities. While previous work has compared principal component analysis (PCA) with KOT methods both directly\cite{mez05, row09, bru16, klu18a, lu20} and indirectly\cite{bra17, red20}, and approaches to do KMD on tensors have been developed \cite{klu18b, gel19}, little work has be done on comparing KMD to TCA \cite{lus19}.  

We considered dynamic mode decomposition (DMD) \cite{row09, sch10, che12, tu14, wil15}, a popular approach for performing KMD, on a data three-tensor, with one dimension being the elements of the state space, one being time, and one being the initial conditions. We proved, in Lemma 2, that when the data is linearly consistent and the Koopman operator, approximated via Exact DMD \cite{tu14}, has a full set of eigenvectors and rank R, the DMD modes are an R--optimal TCA decomposition of the data. Motivated by this, we then formulated a correspondence between the TCA and DMD modes (Eq. 17). We proved, in Lemma 3, that this correspondence was exact, up to scaling and permutation of labels of the modes, when R is equal to the tensor rank of the data and when a certain inequality (Eq. \ref{Kruskal uniqueness}), known to be sufficient for TCA uniqueness, is satisfied. When the modes of the two methods are equivalent, there is a strong implication about the dynamics of the data. Namely, the data comes from distinct sources with single exponential growth/decay and/or oscillatory time dynamics. On a simple example, we showed that modes computed via a numerical implementation of TCA\cite{tensorlab} very nearly matched the modes computed via Exact DMD \cite{tu14} (Fig. \ref{r2_lin_sys}).

\subsection{Limitations}
Our present analysis is limited by the fact that we used Exact DMD to make the connection between TCA and KMD. This requires the approximated Koopman operator to have a full set of eigenvectors and the data to be linearly consistent\cite{tu14}. Because these assumptions are not always met \cite{bag13}, and because there are other DMD and KMD approaches that still hold in such cases, future work will need to look at connecting TCA and KMD in less restricted settings.

In addition, DMD is ultimately limited by the fact that it is an approximation to the true KMD. This approximation can break down when the data's dynamics are nonlinear, which has motivated the development of more sophisticated KOT methods \cite{wil15, kam20, tak17, li17, lus18, ott19, yeu19}. Understanding whether and, if so, how such techniques can be connected to TCA is a major open question, whose answer(s) will greatly extend our work. 

\subsection{TCA vs. DMD}
A large, and growing, body of applied research makes it clear that TCA and DMD are powerful tools, especially when the equations underlying the data's dynamics are unknown. However, because they have largely not been compared to each other \cite{lus19}, it has been unclear whether settings in which one method was used may have benefited from the use of the other method. An important part of this paper is that, by making the connection between the two methods explicit, this kind of comparison can be more easily carried out. Below, we discuss points that have emerged, some of which requires future work to confirm.

First, the number of modes needed to optimally perform DMD will be greater than the number needed for a TCA rank decomposition [i.e. $R > R^*(\pmb{\mathscr{Y}})$]. This is because the DMD modes take an explicit form in their time dependence, whereas the TCA modes do not. Therefore, TCA can be very useful for succinctly describing the data when the sources underlying the data are highly nonlinear. 

Second, work using TCA has emphasized the fact that the modes can be used to interpret and understand the system under study. This is in some contrast to work using DMD, where the emphasis has largely been on prediction. Techniques implemented by those using TCA may therefore be of use to DMD practitioners (see Sec. VIC for discussion of one such technique). Additionally, this implies that, with the correspondence of Eq. 17 in hand, those familiar with TCA should be able to (more) easily interpret DMD results. 

Interestingly (given this last point), DMD may provide advantages over TCA in generating modes that are interpretable. There are several reasons for this. First, depending on the problem, there may be many local minima in the landscape of the objective function used to find the TCA decomposition (as was conjectured to be the case in some instances of the numerical example in Sec. V). This means that each application of TCA may result in different modes, potentially complicating interpretation. The DMD modes, on the other hand, are found via a single set of computations. The resulting decomposition can be unique if the (reduced) singular values of the data matrices are distinct. Thus, even when $R > R^*(\pmb{\mathscr{Y}})$, if the assumptions of the present analysis are met, DMD can provide a unique set of modes that exactly reconstruct the data tensor. This also sidesteps TCA's necessary, and possibly expensive, search for $R^*(\pmb{\mathscr{Y}})$. In addition to extending the uniqueness of TCA, the time dependencies of the DMD modes have exponential and/or oscillatory dynamics. This may make them simpler to interpret in the context of the system of study than the more complicated (nonlinear) modes of the TCA rank decomposition. In support of this, work that modeled the time dependencies of each mode as a sparse linear combination of functions from an over-complete library found success in improving the interpretability of the modes \cite{lus19}. Unlike that approach, DMD does not allow for non-exponential or oscillatory modes, nor does it allow for linear combinations of exponential and/or oscillatory modes [e.g. $\sin(x) + \sin(2x)$]. 

Lastly, DMD may be faster than TCA. This is because DMD requires only matrix operations, whereas TCA requires (many) sequential iterations before convergence, for each of the (many) randomly initialized starting points (but see Battaglino et al. 2018\cite{bat18} and Erichson et al. 2020\cite{eri20} for fast randomized TCA approaches). Indeed, it is known that doing TCA is a non-convex optimization problem \cite{ram20}. In support of this, DMD was found to be faster than PCA when used for prediction \cite{lu20}.

\subsection{Future directions}
By finding connections between TCA and KMD, our work opens up a number of new directions. First, with this correspondence in hand, it should become easier for those familiar with Koopman operator theory (KOT) to communicate with those who work in fields where PCA and TCA are ``gold standard'' tools. These fields include neuroscience and biology, areas where KOT has only recently begun to be applied \cite{bru16, mar20, bal20}. Second, TCA implementations now become a possible tool in KOT's growing numerical toolbox\cite{tu14, mez20}. Most implementations of TCA make use of optimization algorithms \cite{car70, har70, tensorlab, kol09, hon20, ram20}, which are largely absent from the KOT literature (although, relatedly, there has become an increasing interest in using neural networks to perform KMD \cite{tak17, li17, lus18, ott19, yeu19}). Whether these methods are advantageous over the current KMD methods in more general scenarios is an important open question. Additionally, because for many systems it does not make physical sense to have elements of the state space take negative values in their participation in each mode, non-negative TCA has been developed\cite{paa94, lee99, qi16}. This approach constrains all the modes to have positive $\textbf{a}_r$ and has allowed for simpler, more intuitive interpretation of the data. To our knowledge, no similar method exists for KMD. Third, whether using the DMD modes, when their error in reconstructing the data is small, but non-zero, as the starting guess for the TCA optimization algorithm would lead to a faster convergence onto ``good'' TCA modes is another open question that we plan on addressing. In support of this idea, higher order SVD has been found to sometimes be a good starting point for the alternating least squares approach to TCA\cite{kol09}. And fourth, whether and how our results change when the underlying dynamical system is non-autonomous, an area of recent active research in the KOT community \cite{mac18}, will be the direction of future work.

\subsection{KOT and Optimization}
Lastly, we see this work as being part of a growing body of literature that is bringing attention to the fact that dynamical systems theory, and in particular KOT, can be used for problems that have historically relied on optimization theory \cite{dog20, man20, tan20, die20, sah20}. These papers have highlighted the fact that, while optimization theory has its strengths, its de-emphasis on the history of the system (e.g. gradient descent re-computing the gradient anew at each time step) is a considerable cost that KOT avoids. Here we showed that, in certain scenarios, TCA, which has been formulated as a non-convex optimization problem, gives the same (theoretical) decomposition as DMD. 

\begin{acknowledgments}
We would like to thank: Prof. Igor Mezic for the introduction to KOT and the ensuing insightful conversations; Dean Huang for the helpful input and the catching of an error in the equations; Akshunna S. Dogra for the advice on the best framing of the correspondence of Eq. \ref{TCA KMD corr} and the potential use of DMD as a starting point for TCA; and Cory Brown, as well as the members of the Goard Lab at UCSB, for the discussions on neuroscience and KOT. Lastly, we thank the reviewers for their significant feedback on flaws in the original arguments and ways in which the presentation of the material could be improved. The author is in-part supported by a Chancellor Fellowship from UCSB. 
\end{acknowledgments}

\section*{Data Availability}
The data and code used for the example in Sec. V have been made freely publicly available \cite{wtr_git}. 

\nocite{*}

\end{document}